\documentclass[graybox]{svmult}

\usepackage{type1cm}        
\usepackage{makeidx}         
\usepackage{graphicx}        
\usepackage{multicol}        
\usepackage[bottom]{footmisc}

\usepackage{newtxtext}       %
\usepackage{newtxmath}       
\usepackage[normalem]{ulem}

\makeindex             

\newcommand{\nl}{\\[1.1ex]}
\newcommand{\ba}[1]{\left\langle #1 \right\rangle}
\newcommand{\bb}[1]{\left( #1 \right)}
\newcommand{\sbr}[1]{\left[ #1 \right]}
\newcommand{\iO}{\int _{\Omega _K}}

\newcommand{\divv}{\mathrm{div}}

\newcommand{\dO}{\mathrm{d}\Omega _K}
\newcommand{\uu}{\boldsymbol{u}}
\newcommand{\vv}{\boldsymbol{v}}
\newcommand{\nn}{\boldsymbol{n}}
\newcommand{\diop}{\text{d}}

\newcommand{\localL}[1]{\begin{array}{lll}
\mathbb{M}^{(1)} _{\mathbb{A}^{-1}, K_ #1} & {\mathbb{E}^{2,1}}^T \nl
\mathbb{E}^{2,1} & 0
\end{array}}
\newcommand{\localLK}[1]{\begin{array}{lll}
\mathbb{M}^{(1)} _{\mathbb{A}^{-1}, K} & {\mathbb{E}^{2,1}}^T \nl
\mathbb{E}^{2,1} & 0
\end{array}}

\newcommand{\VJ}[1]{\textcolor{black}{#1}}
\newcommand{\VJTwo}[1]{\textcolor{black}{#1}}
\newcommand{\En}{\mathbb{E}_\mathbb{N}}

\begin{document}

\title*{A conservative hybrid method for Darcy flow}
\author{Varun Jain, Yi Zhang, Jo\"{e}l Fisser, Artur Palha and Marc Gerritsma}
\institute{Varun Jain, Yi Zhang, Jo\"{e}l Fisser, Artur Palha and Marc Gerritsma, \at Faculty of Aerospace Engineering, TU Delft, Kluyverweg 1, 2629 HS, Delft,\\
\email{ $\{$V.Jain,Y.Zhang-14,J.M.Fisser,A.PalhaDaSilvaClerigo,M.I.Gerritsma$\}$@tudelft.nl.}}
%
%
\maketitle

\abstract*{Each chapter should be preceded by an abstract (no more than 200 words) that summarizes the content. The abstract will appear \textit{online} at \url{www.SpringerLink.com} and be available with unrestricted access. This allows unregistered users to read the abstract as a teaser for the complete chapter.
Please use the 'starred' version of the \texttt{abstract} command for typesetting the text of the online abstracts (cf. source file of this chapter template \texttt{abstract}) and include them with the source files of your manuscript. Use the plain \texttt{abstract} command if the abstract is also to appear in the printed version of the book.}

\abstract{We present a hybrid mimetic spectral element formulation for Darcy flow.
The discrete representations for 1) conservation of mass, and 2) inter-element continuity, are topological relations \VJTwo{that lead to sparse matrix systems}.
These constraints are independent of the element size and shape, and thus invariant under mesh transformations.
The resultant algebraic system is extremely sparse even for high degree polynomial basis.
Furthermore, the system can be efficiently assembled and solved for each element separately.
}
\section{Introduction}
\label{sec:intro}
Hybrid formulations \cite{2010Boffi,2015Cockburn,2018Zhang} are classical domain decomposition methods which reduce the problem of solving one global system to many small local systems.
The local systems can then be efficiently solved independently of each other in parallel.

In this work we present a hybrid mimetic spectral element formulation to solve Darcy flow.
We follow \cite{2017Jain} \VJTwo{which} render the constraints on divergence of mass flux, the pressure gradient and the inter-element continuity metric free.
The resulting system is extremely sparse and shows a \VJ{reduced growth in condition number as compared to non-hybrid system.}

This document is structured as follows:
In Section \ref{sec:weak_formulation} we define the weak formulation for Darcy flow.
The basis functions are introduced in Section \ref{sec:basis}.
The evaluation of weighted inner product and duality pairings are discussed in Section \ref{sec:inner}.
In Section \ref{sec:disc} we discuss the formulation of discrete algebraic system.
In Section \ref{sec:results} we present results for a test case taken from \cite{2008Herbin}.
%
\section{Darcy flow formulation}
\label{sec:weak_formulation}
For $\Omega \in \mathbb{R}^d$, where $d$ is the dimension of the domain, the governing equations for Darcy flow, are given by,
\begin{equation*}
\left\lbrace
\begin{array}{ll}
\boldsymbol{u} + \mathbb{A}\ \nabla p & = 0 \nl
\nabla \cdot \boldsymbol{u} & = f
\end{array} \right. \quad \text{in}\ \Omega \quad \text{and} \quad \left\lbrace \begin{array}{lll}
\Gamma & = \Gamma _D \cup \Gamma _N \nl
p & = \hat{p} & \text{on}\ \Gamma _D \nl
\boldsymbol{u}\cdot \boldsymbol{n} & = \hat{\boldsymbol{u}}_{\boldsymbol{n}} & \text{on}\ \Gamma _N
\end{array} \right. \;,
\end{equation*}
where, $\boldsymbol{u}$ is the velocity, $p$ is the pressure, $f$ the \VJTwo{prescribed} source term, $\mathbb{A}$ is a $d \times d$ \VJTwo{symmetric} positive definite matrix, $\hat{p}$ and $\hat{\boldsymbol{u}}_{\boldsymbol{n}}$ are the \VJTwo{prescribed} pressure and flux boundary conditions, respectively.
\nl
\emph{\underline{Notations}}

For $f,g \in L^2\bb{\Omega}$, $\bb{f,g}_\Omega$ denotes the usual $L^2$ - inner product. 

For vector-valued function\VJ{s} in $L^2$ we define the weighted inner product by,
\begin{equation}
\label{eq:weighted}
\bb{\boldsymbol{u},\boldsymbol{v}}_{\mathbb{A}^{-1},\Omega} = \int _\Omega \bb{\boldsymbol{u},\mathbb{A}^{-1} \boldsymbol{v}} d\Omega \;,
\end{equation}
\VJTwo{where $\bb{\cdot \;, \cdot}$ denotes the pointwise inner product.}
\\
Duality pairing, denoted by $\ba{\cdot , \cdot}_{\Omega}$, is the outcome of a linear functional on $L^2 \bb{\Omega}$ acting on elements from $L^2\bb{\Omega}$.

Let $\Omega _K$ be a \VJ{disjoint} partitioning of \VJ{$\Omega$ with total number of elements $K$, and $K_i$ is any element in $\Omega _K$, such that, $K_i \in \Omega _K$}.
We define the following broken Sobolev spaces \cite{2016Carsten},
\VJ{$H\bb{\divv ;\Omega _K} = \prod_{K} H\bb{\divv; K_i}$}, and $H^{1/2} \bb{\partial \Omega _K} = \prod _{K} H^{1/2}\bb{\partial K_i}$.
\nl
\emph{\underline{Weak formulation}}

The Lagrange functional for Darcy flow is defined as,
\begin{equation*}
\begin{array}{ll}
\mathcal{L}\bb{\boldsymbol{u}, p, \lambda;f} = & \frac{1}{2}\iO \boldsymbol{u}^T\mathbb{A}^{-1} \boldsymbol{u}\ \dO + \iO p \bb{\nabla \cdot \uu - f} \dO \nl
&+ \int _{\partial \Omega _K \setminus \partial \Omega} \lambda \bb{\uu \cdot \nn}\ \mathrm{d} \Gamma - \int _{\Gamma _D} \hat{p} \bb{\uu \cdot \nn}\ \mathrm{d} \Gamma - \int _{\Gamma _N} \lambda \bb{\hat{\uu} \cdot \nn}\ \mathrm{d}\Gamma
\end{array} \;.
\end{equation*}
The variational problem is then given by: For given \VJ{$f \in L^2\bb{\Omega _K}$}, find $\uu \in H(\divv ;\Omega _K)$, \VJ{$p \in L^2\bb{\Omega _K}$}, \VJ{$\lambda \in H^\frac{1}{2}\bb{\partial \Omega _K}$},  such that,
\begingroup\makeatletter\def\f@size{9.3}\check@mathfonts
\begin{equation}
\label{eq:var1}
\left\lbrace
\begin{array}{lllll}
\VJ{\bb{\boldsymbol{v},\uu}_{\mathbb{A}^{-1},\Omega _K}} + \ba{\nabla \cdot \vv, p}_{\Omega _K} + \ba{\bb{\vv \cdot \boldsymbol{n}} , \lambda}_{\partial \Omega_K \setminus \partial \Omega} & = \ba{\vv \cdot \nn, \hat{p}}_{\Gamma _D} & \quad \forall\ \vv \in H(\divv ;\Omega _K) \nl
\ba{q, \nabla \cdot \uu}_{\Omega _K} & = \ba{q,f}_{\Omega _K} & \quad \forall\ q \in L^2\bb{\Omega _K} \nl
\ba{\mu , \bb{\uu\cdot \nn} }_{\partial \Omega_K \setminus \partial \Omega} & = \ba{\mu, \hat{\uu} \cdot \nn}_{\Gamma _N} & \quad \forall\ \mu \in H^\frac{1}{2}\bb{\partial \Omega _K}
\end{array} \right. \;.
\end{equation}
\endgroup

\section{Basis functions}
\label{sec:basis}
\emph{\underline{Primal and dual nodal degrees of freedom}}

Let $\xi_j$, $ j = 0, 1, ... ,N$, be the $N+1$ Gauss-Lobatto-Legendre (GLL) points in $I \in \left[-1,1\right]$.
The Lagrange polynomials $h_i(\xi)$ through $\xi_j$, of degree $N$, given by,
\begin{equation*}
h_i \bb{\xi}= \frac{\bb{\xi ^2 -1}L'_N\bb{\xi}}{N\bb{N+1}L_N\bb{\xi _i}\bb{\xi -\xi_i}} \;,
\end{equation*}
form the 1D primal nodal polynomials which satisfy, $h_i(\xi_j) = \delta_{ij}$.

Let $a^h$ and $b^h$ be two polynomials expanded in terms of $h_i\bb{\xi}$.
The $L^2$ - inner product is then given by,
\begin{equation*}
\bb{a^h,b^h} _I = \vec{a}^T \mathbb{M}^{(0)} \vec{b} \;, \quad \text{where} \quad \mathbb{M}^{(0)} _{i,j} = \int _{-1}^1 {h_i} (\xi)\ {h_j} (\xi)\ \mathrm{d}\xi \;,
\end{equation*}
and, $\vec{a} = \left[ \mathsf{a}_0\ \mathrm{a}_1\ \ldots\  \mathsf{a}_{N}\right]$ and $\vec{b} = \left[ \mathsf{b}_0\ \mathsf{b}_1\ \ldots\  \mathsf{b}_{N}\right]$ are the nodal degrees of freedom.
We define the algebraic \emph{dual} degrees of freedom, $\widetilde{\vec{a}}$, such that the duality pairing is simply the vector dot product between primal and dual degrees of freedom,
\begin{equation*}
\left\langle a^h, b^h \right\rangle _{I} = \widetilde{\vec{a}}^T \vec{b} := \vec{a}^T \mathbb{M}^{(0)} \vec{b} \Rightarrow \widetilde{\vec{a}} = \mathbb{M}^{(0)} \vec{a} \;.
\end{equation*}
Thus, the dual degrees of freedom are linear functionals of primal degrees of freedom.
\nl
\emph{\underline{Primal and dual edge degrees of freedom}}

The edge polynomials, for the $N$ edges between $N+1$ GLL points, of polynomial degree $N-1$, are defined as \cite{2011Gerritsma},
\begin{equation*}
    e_j(\xi) = -\sum_{k=1}^{j-1}\frac{\diop h_k}{\diop \xi}(\xi) \;, \quad \text{such that} \quad \int _{\xi _{j-1}} ^{\xi _j} e_i (\xi) = \delta _{ij} \;.
\end{equation*}
Let $p^h$ and $q^h$ be two polynomials expanded in edge basis functions.
%
The inner product in $L^2$ space is given by,
\begin{equation*}
\bb{p^h,q^h} _{I} = \vec{p}^T \mathbb{M}^{(1)} \vec{q} \;, \quad \text{where} \quad \mathbb{M}^{(1)} _{i,j} = \int _{-1}^1 {e_i} (\xi)\ e_j (\xi)\ \mathrm{d}\xi \;,
\end{equation*}
and, $\vec{p} = \left[ \mathsf{p}_1\ \mathsf{p}_2\ \ldots\  \mathsf{p}_{N}\right]$ and $\vec{q} = \left[ \mathsf{q}_1\ \mathsf{q}_2\ \ldots\  \mathsf{q}_{N}\right]$ are the edge degrees of freedom.
As before, we define the \emph{dual} degrees of freedom such that,
\begin{equation*}
\left\langle p^h, q^h \right\rangle _{I} = \widetilde{\vec{p}}^T \vec{q} := \vec{p}^T \mathbb{M}^{(1)} \vec{q} \Rightarrow \widetilde{\vec{p}} = \mathbb{M}^{(1)} \vec{p} \;.
\end{equation*}
\VJ{A similar construction can be used for dual degrees of freedom in higher dimension. For construction of the dual degrees of freedom in $2D$ see \cite{2017Jain} and for $3D$ see \cite{2018Zhang2}.}
\nl
\emph{\underline{Differentiation of nodal polynomial representation}}

Let $a^h\bb{\xi}$ be expanded in Lagrange polynomial, then
\begin{equation} \label{eq:deri}
\frac{d}{d\xi}a^h\bb{\xi} = \frac{d}{d\xi}\sum _{i=0}^N \mathsf{a}_i h_i \bb{\xi} =\sum _{i=1}^N \bb{\mathsf{a}_i - \mathsf{a}_{i-1}} e_i \bb{\xi} \;.
\end{equation}
Therefore, taking \VJTwo{the} derivative of \VJTwo{a} polynomial involves two steps :
First, \VJTwo{take the} difference of degrees of freedom; and second, change of basis from nodal to edge \cite{2011Gerritsma}.
\section{Discrete inner product and duality pairing} \label{sec:inner}
For 2D domains, the higher dimensional primal basis are constructed \VJTwo{using the} tensor product of the 1D basis.

For the weak formulation in (\ref{eq:var1}) we expand the velocity $\uu ^h$ in primal edge basis as,
\begin{equation}
\label{eq:exp}
\uu^h\bb{\xi, \eta} = \sum _{i=0}^N \sum _{j=1}^N \mathsf{u_x} _{i,j}\ h_i(\xi)\ e_j(\eta)\ \vec{\hat{\imath}} + \sum _{i=1}^N \sum _{j=0}^N \mathsf{u_y} _{i,j}\ e_i(\xi)\ h_j(\eta)\ \vec{\hat{\jmath}} \;.
\end{equation} \nl
\emph{\underline{Weighted inner product}}

Using (\ref{eq:weighted}) and the expansions in (\ref{eq:exp}), the weighted inner product is evaluated as,
\begin{equation*}
\bb{\vv^h,\uu ^h}_{\mathbb{A}^{-1}, \Omega _K} = \sum _{K}\vec{v}_{K_i}^T\ \mathbb{M}^{(1)} _{\mathbb{A}^{-1},{K_i}}\ \vec{u}_{K_i} \;,
\end{equation*}
where, $\vec{u}_{K_i}$ are the degrees of freedom in element ${K_i}$, and
\begin{equation*}
\mathbb{M}^{(1)} _{\mathbb{A}^{-1},{K_i}} = \int _{K_i} 
\bb{\begin{array}{l}
h_i(\VJ{x})\ e_j(\VJ{y}) \nl
e_i(\VJ{x})\ h_j(\VJ{y})
\end{array}}
\times \mathbb{A}^{-1}\bb{\VJ{x},\VJ{y}}
\bb{\begin{array}{l}
h_i(\VJ{x})\ e_j(\VJ{y}) \nl
e_i(\VJ{x})\ h_j(\VJ{y})
\end{array}}\
\mathrm{d} {K_i} \;.
\end{equation*}
For mapping of elements please refer to \cite{2017Gerritsma}. \nl
\emph{\underline{Divergence of velocity}}

Divergence of velocity, $\nabla \cdot \uu^h$, is evaluated using (\ref{eq:deri}), but now for 2D,
\begin{equation} \label{eq:disc_divu}
\begin{array}{ll}
\nabla \cdot \uu ^h & = \frac{\partial}{\partial x} \sum _{i=0}^N \sum _{j=1}^N \mathsf{u_x} _{i,j}\ h_i(x) e_j(y) + \frac{\partial}{\partial y} \sum _{i=1}^N \sum _{j=0}^N \mathsf{u_y} _{i,j}\ e_i(x) h_j(y) \nl
& = \sum _{i,j =1} ^N \bb{\mathsf{u_x}_{i,j} - \mathsf{u_x}_{i-1,j} + \mathsf{u_y}_{i,j} - \mathsf{u_y}_{i,j-1}} e_i\bb{x} e_j \bb{y}
\end{array} \;.
\end{equation}
The pressure is expanded in the dual basis $\widetilde{e_i\bb{\xi}e_j\bb{\eta}}$.
These basis are dual to the basis in which $\nabla \cdot \uu ^h$ is expanded in (\ref{eq:disc_divu}).
Therefore the \VJTwo{weak} constraint on divergence of velocity is a duality pairing evaluated as,
\begin{equation*}
\VJ{\ba{q^h, \nabla \cdot \uu^h}_ {\Omega_K} } = \sum _{K} \vec{q}_{K_i}^T\ \mathbb{E}^{2,1}\ \vec{u}_{K_i} \;,
\end{equation*}
where $\mathbb{E}^{2,1}$ represents the discrete divergence operator.
It is an incidence matrix that is metric-free and topological, and remains the same for each element \VJ{in} $\Omega _K$.
For an extensive discussion on the incidence matrix, see \cite{2017Gerritsma}.
For an element of degree $N=3$,
\begingroup\makeatletter\def\f@size{6}\check@mathfonts
\begin{equation*}
\mathbb{E}^{2,1} = \left[ \begin{array}{cccccccccccccccccccccccc}
-1 & 0 & 0 & 1 & 0 & 0 & 0 & 0 & 0 & 0 & 0 & 0 & -1 & 0 & 0 & 1 & 0 & 0 & 0 & 0 & 0 & 0 & 0 & 0 \\
0  & 0 & 0 & -1 & 0 & 0 & 1 & 0 & 0 & 0 & 0 & 0 & 0 & -1 & 0 & 0 & 1 & 0 & 0 & 0 & 0 & 0 & 0 & 0 \\
0  & 0 & 0 & 0  & 0 & 0 & -1 & 0 & 0 & 1 & 0 & 0 & 0 & 0 & -1 & 0 & 0 & 1 & 0 & 0 & 0 & 0 & 0 & 0 \\
0 & -1 & 0 & 0 & 1 & 0 & 0 & 0 & 0 & 0 & 0 & 0 & 0 & 0 & 0 & -1 & 0 & 0 & 1 & 0 & 0 & 0 & 0 & 0 \\
0 & 0  & 0 & 0 & -1 & 0 & 0 & 1 & 0 & 0 & 0 & 0 & 0 & 0 & 0 & 0 & -1 & 0 & 0 & 1 & 0 & 0 & 0 & 0 \\
0 & 0 & 0 & 0 & 0 & 0 & 0 & -1 & 0 & 0 & 1 & 0 & 0 & 0 & 0 & 0 & 0 & -1 & 0 & 0 & 1 & 0 & 0 & 0 \\
0 & 0 & -1 & 0 & 0 & 1 & 0 & 0 & 0 & 0 & 0 & 0 & 0 & 0 & 0 & 0 & 0 & 0 & -1 & 0 & 0 & 1 & 0 & 0 \\
0 & 0 & 0 & 0 & 0 & -1 & 0 & 0 & 1 & 0 & 0 & 0 & 0 & 0 & 0 & 0 & 0 & 0 & 0 & -1 & 0 & 0 & 1 & 0 \\
0 & 0 & 0 & 0 & 0 & 0 & 0 & 0 & -1 & 0 & 0 & 1 & 0 & 0 & 0 & 0 & 0 & 0 & 0 & 0 & -1 & 0 & 0 & 1
\end{array} \right] \;.
\end{equation*}
\endgroup
\nl
\emph{\underline{Connectivity matrix}}

\VJ{The connectivty matrix ensures continuity of the velocity across the elements.
$\vec{\lambda}$ is the interface variable between the elements that acts as Lagrange multiplier that imposes the constraint given by,
\begin{equation*} \label{eq:en}
    \ba{\mu ^h,\uu^h \cdot \nn}_{\partial \Omega _K \setminus \partial \Omega} = \sum _K \vec{\mu}_{K_i}^T\ \mathbb{N}\ \vec{u}_{K_i}  = \vec{\mu}^T\ \mathbb{E}_\mathbb{N}\ \vec{u}\ \;,
\end{equation*}
\VJTwo{where $\mathbb{N}$ is the discrete trace operator.
It is a sparse matrix that consists of $1$, $-1$ and $0$ only.}
\VJTwo{For construction of $\mathbb{N}$ please refer to \cite{2018Gerritsma}}.
\VJTwo{$\mathbb{E}_\mathbb{N}$ is the assembled $\mathbb{N}$ for all the elements.}
For discretization, $K = 2 \times 2$, $N=2$, $\mathbb{E}_\mathbb{N}$ is shown in (\ref{eq:en2}).
The matrix size of $\mathbb{E}_\mathbb{N}$ is $8 \times 64$, but it has only 16 non-zero entities.
It is an extremely sparse matrix that is metric-free and the location of +/-1 valued entries depend only on the connection between different elements.}

\begin{equation} \label{eq:en2}
\mathbb{E}_\mathbb{N} = 
\begingroup\makeatletter\def\f@size{1}\check@mathfonts
\left[ \begin{array}{lllllllllllllllllllllllllllllllllllllllllllllllllllllllllllllllllll}
0&	0&	0&	0&	\scalebox{0.5}{1}&	0&	0&	0&	0&	0&	0&	0&	0&	0&	0&	0&	\scalebox{0.5}{-1}&	0&	0&	0&	0&	0&	0&	0&	0&	0&	0&	0&	0&	0&	0&	0&	0&	0&	0&	0&	0&	0&	0&	0&	0&	0&	0&	0&	0&	0&	0&	0&	0&	0&	0&	0&	0&	0&	0&	0&	0&	0&	0&	0&	0&	0&	0&	0 \\
0&	0&	0&	0&	0&	\scalebox{0.5}{1}&	0&	0&	0&	0&	0&	0&	0&	0&	0&	0&	0&	\scalebox{0.5}{-1}&	0&	0&	0&	0&	0&	0&	0&	0&	0&	0&	0&	0&	0&	0&	0&	0&	0&	0&	0&	0&	0&	0&	0&	0&	0&	0&	0&	0&	0&	0&	0&	0&	0&	0&	0&	0&	0&	0&	0&	0&	0&	0&	0&	0&	0&	0 \\
0&	0&	0&	0&	0&	0&	0&	0&	0&	0&	0&	0&	0&	0&	0&	0&	0&	0&	0&	0&	0&	0&	0&	0&	0&	0&	0&	0&	0&	0&	0&	0&	0&	0&	0&	0&	\scalebox{0.5}{1}&	0&	0&	0&	0&	0&	0&	0&	0&	0&	0&	0&	\scalebox{0.5}{-1}&	0&	0&	0&	0&	0&	0&	0&	0&	0&	0&	0&	0&	0&	0&	0 \\
0&	0&	0&	0&	0&	0&	0&	0&	0&	0&	0&	0&	0&	0&	0&	0&	0&	0&	0&	0&	0&	0&	0&	0&	0&	0&	0&	0&	0&	0&	0&	0&	0&	0&	0&	0&	0&	\scalebox{0.5}{1}&	0&	0&	0&	0&	0&	0&	0&	0&	0&	0&	0&	\scalebox{0.5}{-1}&	0&	0&	0&	0&	0&	0&	0&	0&	0&	0&	0&	0&	0&	0 \\
0&	0&	0&	0&	0&	0&	0&	0&	0&	0&	\scalebox{0.5}{1}&	0&	0&	0&	0&	0&	0&	0&	0&	0&	0&	0&	0&	0&	0&	0&	0&	0&	0&	0&	0&	0&	0&	0&	0&	0&	0&	0&	\scalebox{0.5}{-1}&	0&	0&	0&	0&	0&	0&	0&	0&	0&	0&  	0&	0&	0&	0&	0&	0&	0&	0&	0&	0&	0&	0&	0&	0&	0 \\
0&	0&	0&	0&	0&	0&	0&	0&	0&	0&	0&	\scalebox{0.5}{1}&	0&	0&	0&	0&	0&	0&	0&	0&	0&	0&	0&	0&	0&	0&	0&	0&	0&	0&	0&	0&	0&	0&	0&	0&	0&	0&	0&	\scalebox{0.5}{-1}&	0&	0&	0&	0&	0&	0&	0&	0&	0&	0&	0&	0&	0&	0&	0&	0&	0&	0&	0&	0&	0&	0&	0&	0 \\
0&	0&	0&	0&	0&	0&	0&	0&	0&	0&	0&	0&	0&	0&	0&	0&	0&	0&	0&	0&	0&	0&	0&	0&	0&	0&	\scalebox{0.5}{1}&	0&	0&	0&	0&	0&	0&	0&	0&	0&	0&	0&	0&	0&	0&	0&	0&	0&	0&	0&	0&	0&	0&	0&	0&	0&	0&	0&	\scalebox{0.5}{-1}&	0&	0&	0&	0&	0&	0&	0&	0&	0 \\
0&	0&	0&	0&	0&	0&	0&	0&	0&	0&	0&	0&	0&	0&	0&	0&	0&	0&	0&	0&	0&	0&	0&	0&	0&	0&	0&	\scalebox{0.5}{1}&	0&	0&	0&	0&	0&	0&	0&	0&	0&	0&	0&	0&	0&	0&	0&	0&	0&	0&	0&	0&	0&	0&	0&	0&	0&	0&	0&	\scalebox{0.5}{-1}&	0&	0&	0&	0&	0&	0&	0&	0
\end{array} \right] \;.
\endgroup
\end{equation}

\section{Discrete formulation}
\label{sec:disc}
Using the weighted inner product and duality pairings discussed in Section~\ref{sec:inner}, we can write the discrete form of weak formulation in (\ref{eq:var1}) as,
%
\begin{equation}
\label{eq:sys}
\left[
\begin{array}{cc}
\mathbb{A} & {\mathbb{E}_\mathbb{N}}^T \nl
\mathbb{E}_\mathbb{N} & 0
\end{array}
\right]
\left[\begin{array}{l}
\vec{X} \nl
\vec{\lambda}
\end{array}
\right] = \sbr{\begin{array}{l}
\vec{F} \nl
0
\end{array}} \;,
\end{equation}
where, $\mathbb{A}$ is an invertible block diagonal matrix given by,
\begin{equation} \label{eq:A}
\mathbb{A} = \left[ \begin{array}{llllllllllllllllllllll}
\localL{1} \nl
& & \localL{2} \nl
& & & & \ddots & \ddots \nl
& & & & \ddots & \ddots \nl
& & & & & & \localLK{K} 
\end{array} \right] \;,
\end{equation}
$\mathbb{E}_\mathbb{N}$ is as given in (\ref{eq:en2}), $ \vec{X} =  \sum _{K} \sbr{\begin{array}{c}
\vec{u} \nl
\vec{p}
\end{array}}_{K_i}$, and $\vec{F} = \sum _{K} \sbr{\begin{array}{c}
\vec{\hat{p}} \nl
\vec{f}
\end{array}}_{K_i} $, where $\vec{f}$ are the expansion coefficients of \VJ{$f^h\bb{x,y} = \sum _{i,j}^N \mathsf{f}_{ij}\ e_i\bb{x} e_j\bb{y}$}. 

In (\ref{eq:A}), the mass matrix $\mathbb{M}^{(1)} _{\mathbb{A}^{-1},{K_i}}$ is the only dense matrix and also the only component that changes with each local element, ${K_i}$.
$\mathbb{E}_\mathbb{N}$ is a sparse incidence matrix for the global system and $\mathbb{E}^{2,1}$ is a sparse incidence matrix for the local systems that remain the same for each element.

Using the Schur complement method, the global system (\ref{eq:sys}) can be reduced to solve for $\vec{\lambda}$, \cite{2010Boffi},
\begin{equation} \label{eq:lam}
\vec{\lambda} = \bb{\En \mathbb{A}^{-1}\ {\En} ^T}^{-1} \cdot \bb{\En \mathbb{A}^{-1} \vec{F}} \;.
\end{equation}
\VJ{To evaluate $\vec{\lambda}$ in (\ref{eq:lam}) we need $\mathbb{A}^{-1}$ that can be calculated efficiently by taking inverse of each block of $\mathbb{A}$ separately.
This part can be easily parallelized.}
Once the $\vec{\lambda}$ is determined the solution in each element, ${K_i}$, can be evaluated by solving for,
\begin{equation} \label{eq:eval}
\left[\begin{array}{l}
\vec{u}\nl
\vec{p}
\end{array}
\right]_{K_i} = \left[
\begin{array}{ccc}
\mathbb{M}^{(1)} _{\mathbb{A}^{-1}} & {\mathbb{E}^{2,1}}^T \nl
\mathbb{E}^{2,1} & 0
\end{array}
\right]_{K_i} ^{-1} 
\left[ \begin{array}{l}
\vec{\lambda} \nl
\vec{f}
\end{array} \right]_{K_i} \;.
\end{equation}
\VJ{Here the inverse of the local block in the RHS is already evaluated during (\ref{eq:lam}).}
As the local systems are independent of each other (\ref{eq:eval}) can also be evaluated separately for each element \VJ{and easily parallelized.}

\VJ{The system (\ref{eq:lam}) solves for interface degrees of freedom between the elements and will always be smaller than the full global system.
For a comparison of the size of $\lambda$ system with full system see Table 1 for 2D systems, and Table 2 for 3D systems.
In Table 1 \& 2 (left) we see that, for constant $K$, increasing the order of polynomial basis the growth in size of $\lambda$ system is less than the growth in size of full system.
Thus, hybrid formulations are beneficial for high order methods, in 2D and in 3D, where local degrees of freedom of an element are much higher than interface degrees of freedom.}

In Table 1 \& 2 (right) we see that, for constant $N$, the $\lambda$ system is certainly smaller than the full system, although the growth rate in size of $\lambda$ and full systems does not change significantly.
\begin{table}
\label{tab:ndof}
\center
\caption{For 2D. Left: Number of total unknowns as a function of $N$, for $K=3 \times 3$.
Right: Number of unknowns as a function of $K$, for $N=3$.}
\begin{tabular}{|c|c|c|c| }
\hline
 $N$ & Full system & $\lambda$ only & $\lambda$ / Full \\
 \hline
 $5 $ & $825$   & $60$ 	& 0.07 \\  
 $10$ & $3000$  & $120$ 	& 0.04 \\
 $15$ & $6525$  & $180$ 	& 0.03 \\  
 $20$ & $11400$ & $240$ 	& 0.02 \\  
 $25$ & $17625$ & $300$ 	& 0.02 \\  
 \hline
\end{tabular} \quad \quad
\begin{tabular}{|c|c|c|c| }
\hline
 $K$ & Full system & $\lambda$ only & $\lambda$ / Full \\
 \hline
 $400 $  & $15480$  & $2280$ 	& 0.15\\  
 $1600$  & $62160$  & $9360$ 	& 0.15\\
 $3600$  & $140040$ & $21240$ 	& 0.15\\  
 $6400$  & $249120$ & $37920$ 	& 0.15\\  
 $10000$ & $389400$ & $59400$ 	& 0.15\\  
 \hline
\end{tabular}
\end{table}
\vspace{-1cm}
\begin{table}
\label{tab:ndofb}
\center
\caption{For 3D. Left: Number of total unknowns as a function of $N$, for $K=3 \times 3 \times 3$.
Right: Number of unknowns as a function of $K$, for $N=3$.}
\begin{tabular}{|c|c|c|c| }
\hline
 $N$ & Full system & $\lambda$ only & $\lambda$ / Full \\
 \hline 
 $5 $ & $ 16875 $ & $ 1350 $ 	& $ 0.08 $ \\  
 $10$ & $ 121500 $ & $ 5400 $ 	& $ 0.04 $ \\
 $15$ & $ 394875 $ & $ 12150 $ 	& $ 0.03 $ \\  
 $20$ & $ 918000 $ & $ 21600 $ 	& $ 0.02 $ \\  
 $25$ & $ 1771875 $ & $ 33750 $ 	& $ 0.02 $ \\  
 \hline
\end{tabular} \quad \quad
\begin{tabular}{|c|c|c|c| }
\hline
 $K$ & Full system & $\lambda$ only & $\lambda$ / Full \\
 \hline
 $8000 $  & $ 1285200 $ & $ 205200 $ 	& $ 0.16 $\\  
 $64000$  & $ 10324800 $ & $ 1684800 $ 	& $ 0.16 $\\
 $216000$  & $ 34894800 $ & $ 5734800 $ 	& $ 0.16 $\\  
 $512000$  & $ 82771200 $ & $ 13651200 $ 	& $ 0.16 $\\  
 $1000000$ & $ 161730000 $ & $ 26730000 $ 	& $ 0.17 $\\  
 \hline
\end{tabular}
\end{table}
\vspace{-1cm}
\section{Results}
\label{sec:results}
In this section we present the results for a test problem from \cite{2008Herbin} by solving system (\ref{eq:sys}).
The domain of test problem is, $\Omega \in \left[0,1 \right]^2$. 
The source term is defined as,
\begin{equation*}
f_{ex} = \nabla \cdot (-\mathbb{A}\nabla p _{ex}) \;, \quad \text{where} \;,
\end{equation*}
\begin{equation*}
\begin{array}{ll}
\mathbb{A} & = \frac{1}{x^2 + y^2 + \alpha} \bb{\begin{array}{ll}
10^{-3}x^2 + y^2 + \alpha & \quad \bb{10^{-3}-1}xy \nl
\bb{10^{-3}-1}xy & \quad x^2 + 10^{-3}y^2 + \alpha
\end{array}} \;; \quad \quad \alpha = 0.1 \nl
p _{ex} & = \mathrm{sin}\bb{2 \pi x} \mathrm{sin}\bb{2 \pi y} 
\end{array} \;,
\end{equation*}
and Dirichlet boundary conditions are imposed along the entire boundary, $\Gamma _\mathrm{D}= \Gamma$ and $\Gamma _\mathrm{N} = \emptyset $.
\begin{figure}
\center
\includegraphics[scale=0.4]{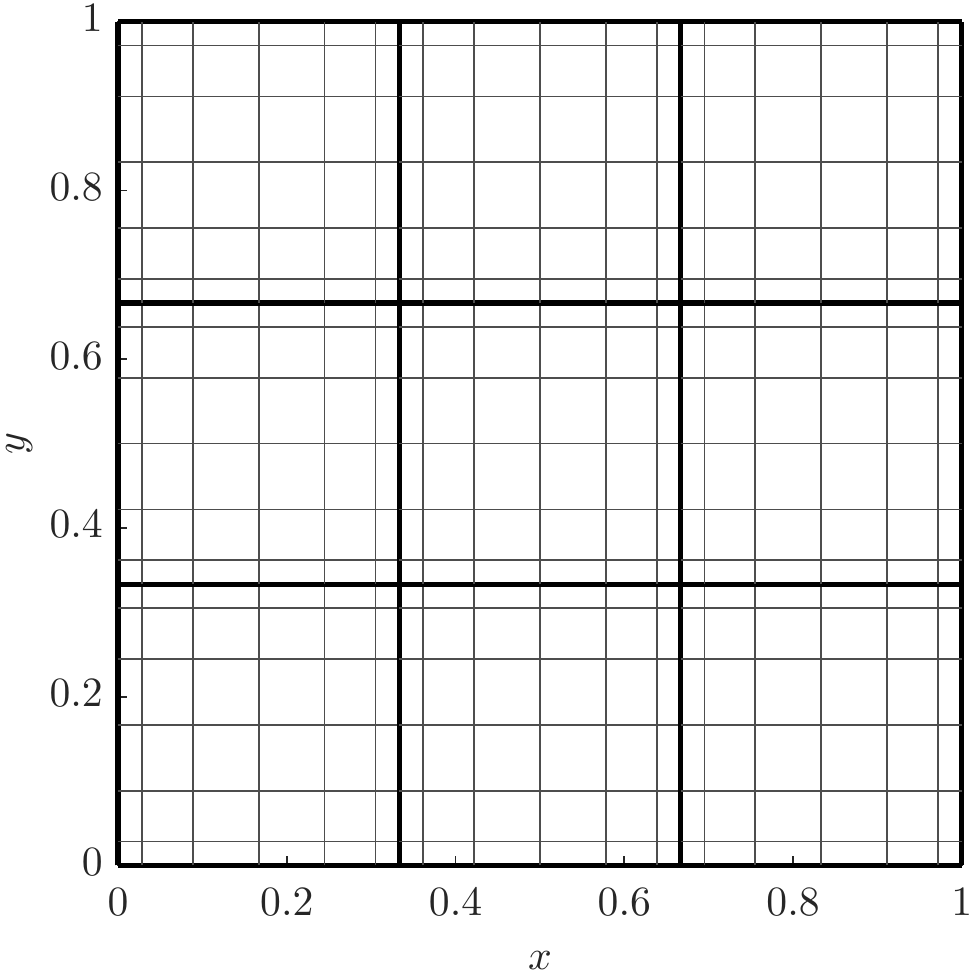}
\includegraphics[scale=0.4]{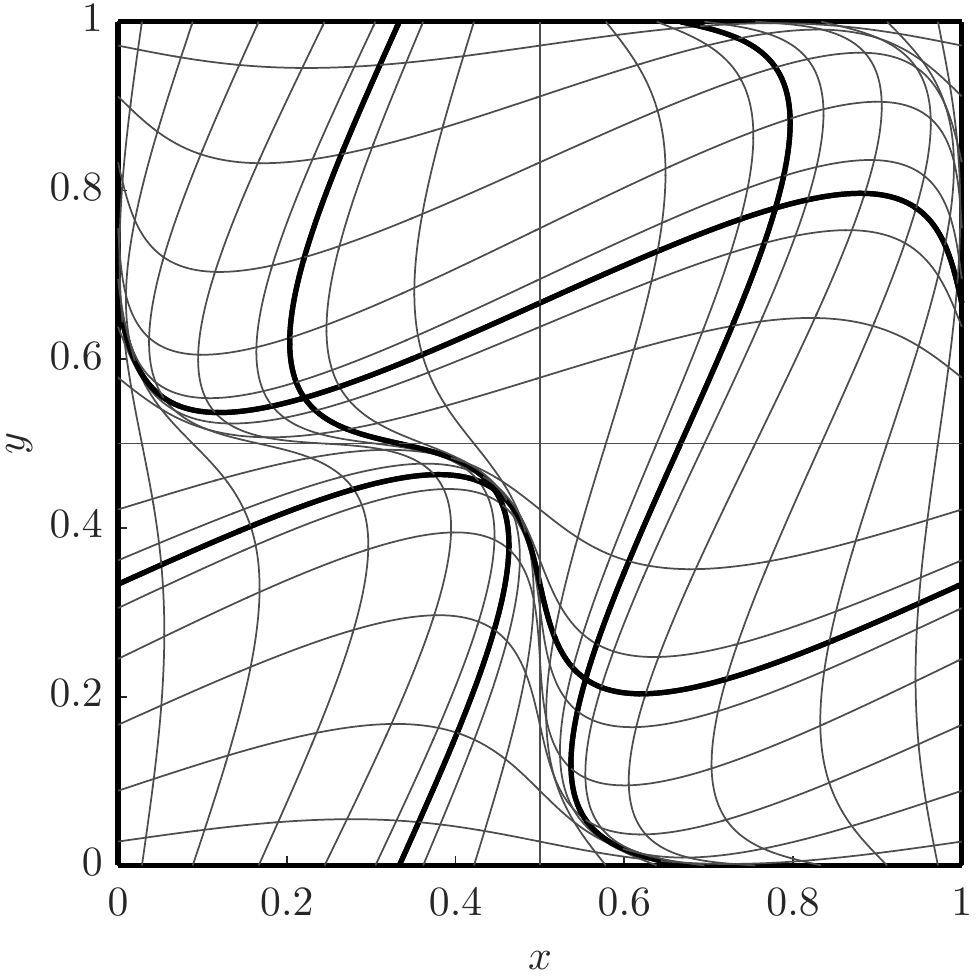}
\caption{Mesh configuration: $K = 3 \times 3$, $N=6$, Left: orthogonal, Right: curved.}
\label{fig:mesh}
\end{figure}
We solve this problem on an orthogonal and a highly curved mesh, see Fig. \ref{fig:mesh}.
\begin{figure}
\center
\includegraphics[scale=0.3]{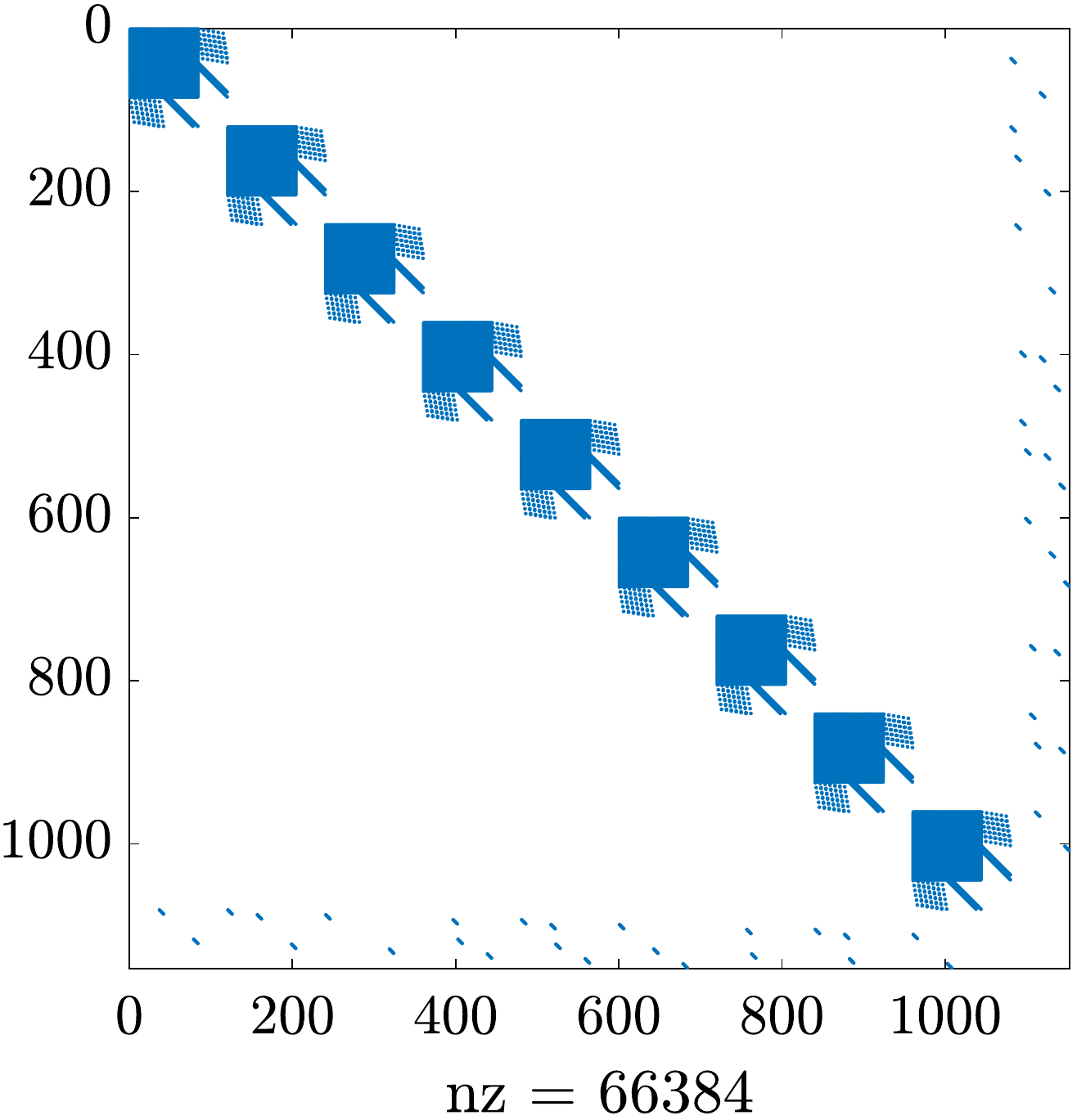} \quad
\includegraphics[scale=0.3]{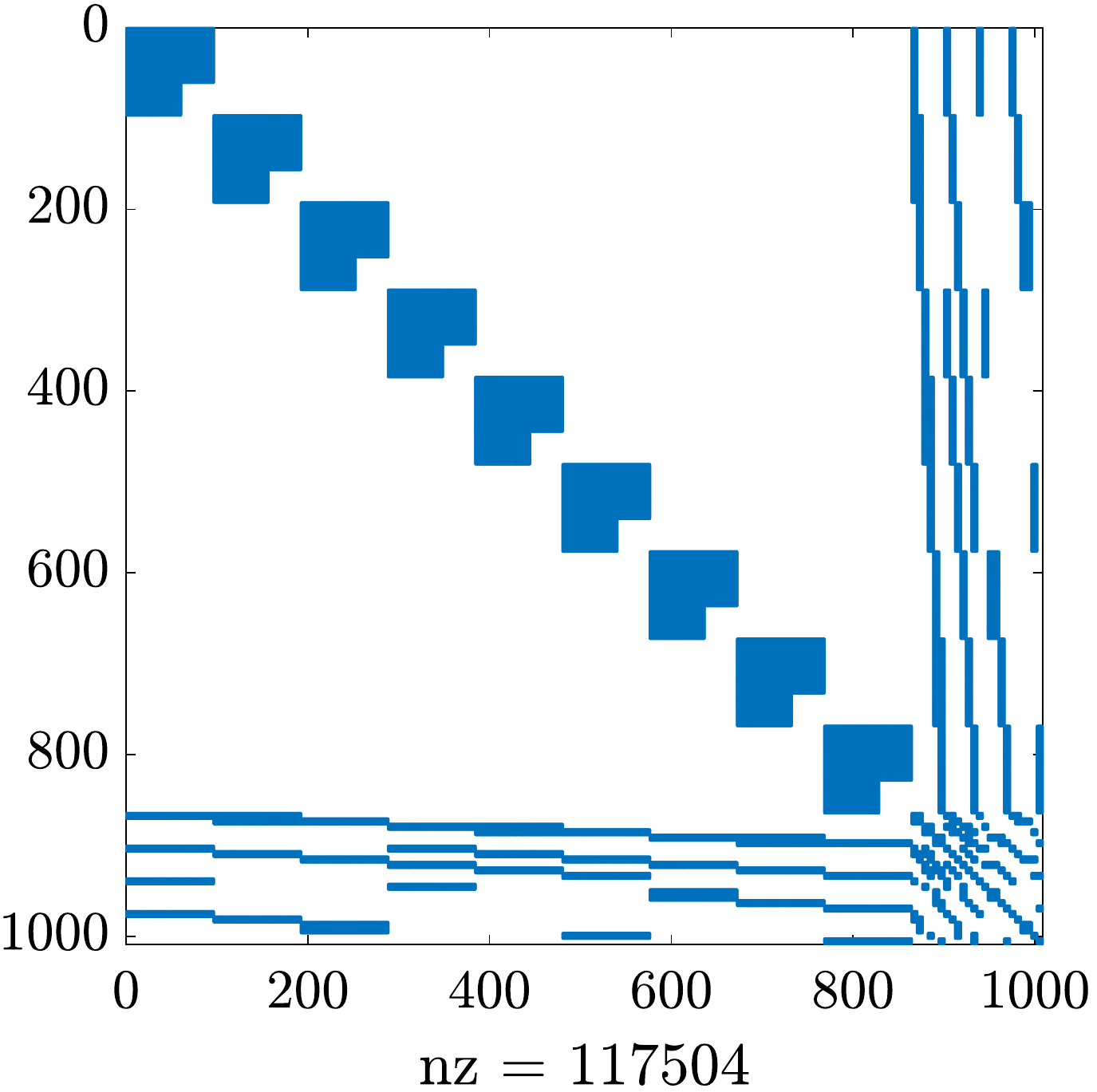}
\caption{Sparsity plots $K = 3 \times 3$, $N=6$. Left: Hybrid elements method. Right: Continuous element method.}
\label{fig:spy}
\end{figure}

The same problem was earlier addressed by authors in \cite{2017Gerritsma}, but for a method with continuous elements and \emph{primal} basis only.
For the configuration $K=3 \times 3$, $N=6$, we compare the sparsity structure of the two approaches in Fig \ref{fig:spy}.
On left we see the hybrid formulation, and on the right we see the continuous elements formulation \VJ{\cite{2017Gerritsma}}.
The number of non zero entities are almost half in the hybrid formulation, $66384$, as compared to the continuous element formulation, $117504$.
Here, the sparsity is due to use of algebraic dual degrees of freedom and is not because of hybridization of the scheme.
%
\begin{figure}
\includegraphics[scale=0.3]{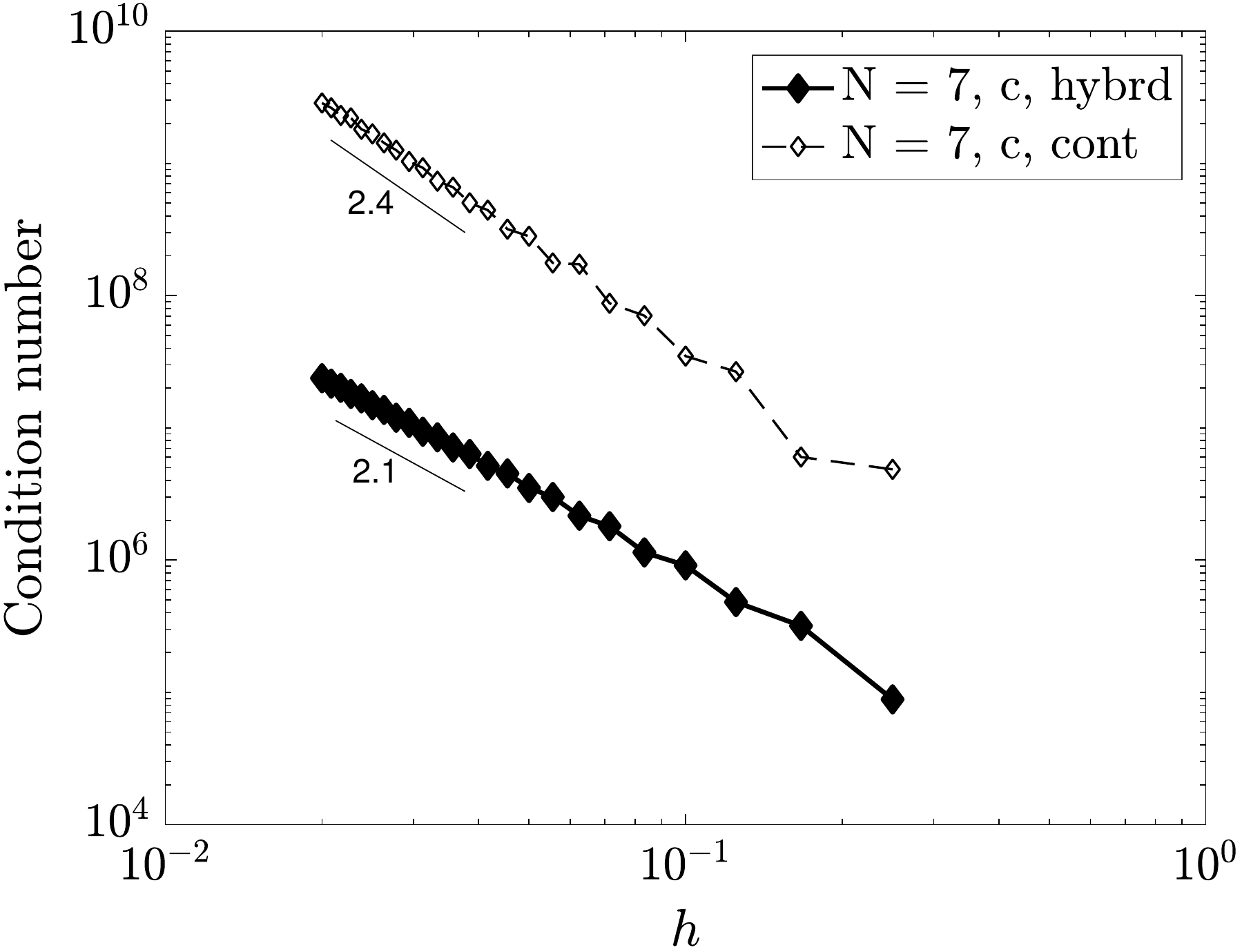}
\includegraphics[scale=0.3]{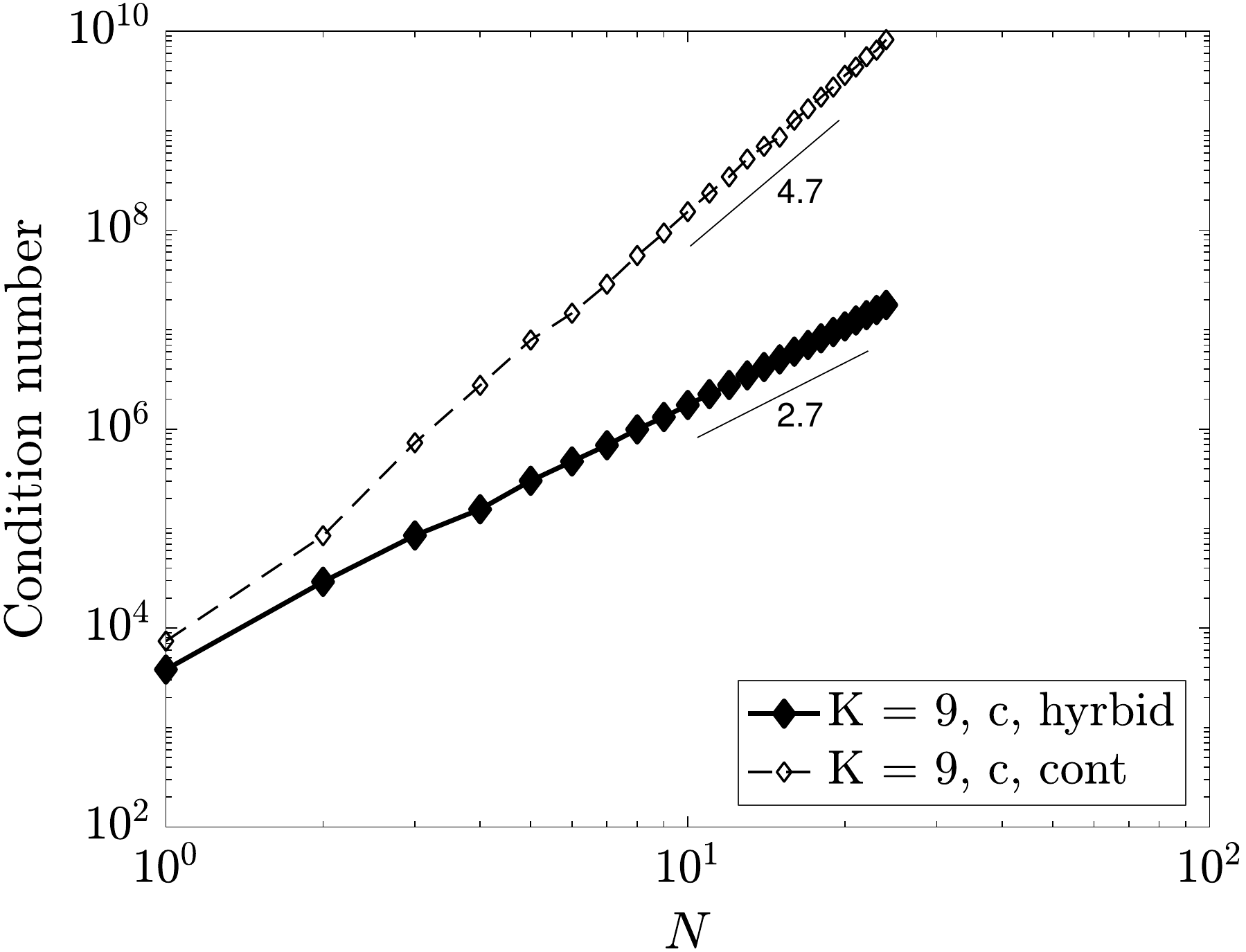}
\caption{Growth in condition number for hybrid elements in dark line, and continuous elements in dotted line. Left: $h$-refinement; Right: $N$-refinement.}
\label{fig:cond}
\end{figure}

In Fig. \ref{fig:cond}, on the left we compare the growth in condition number, for the $\vec{\lambda}$ only system with continuous element system, \VJ{for $N=7$ on the curved mesh,} with increasing number of elements, $K$.
We observe similar growth rates for hybrid and continuous formulation, however the condition number for continuous elements formulation is \VJ{almost $\mathcal{O}\bb{10^2}$} higher.
On the right we see the growth in condition number with increasing polynomial degree \VJ{for $K=9 \times 9$ on the curved mesh}.
A suppressed growth rate in condition number for hybrid formulation is observed.
Thus hybrid formulations are beneficial for high order methods.
\begin{figure}
\includegraphics[scale=0.3]{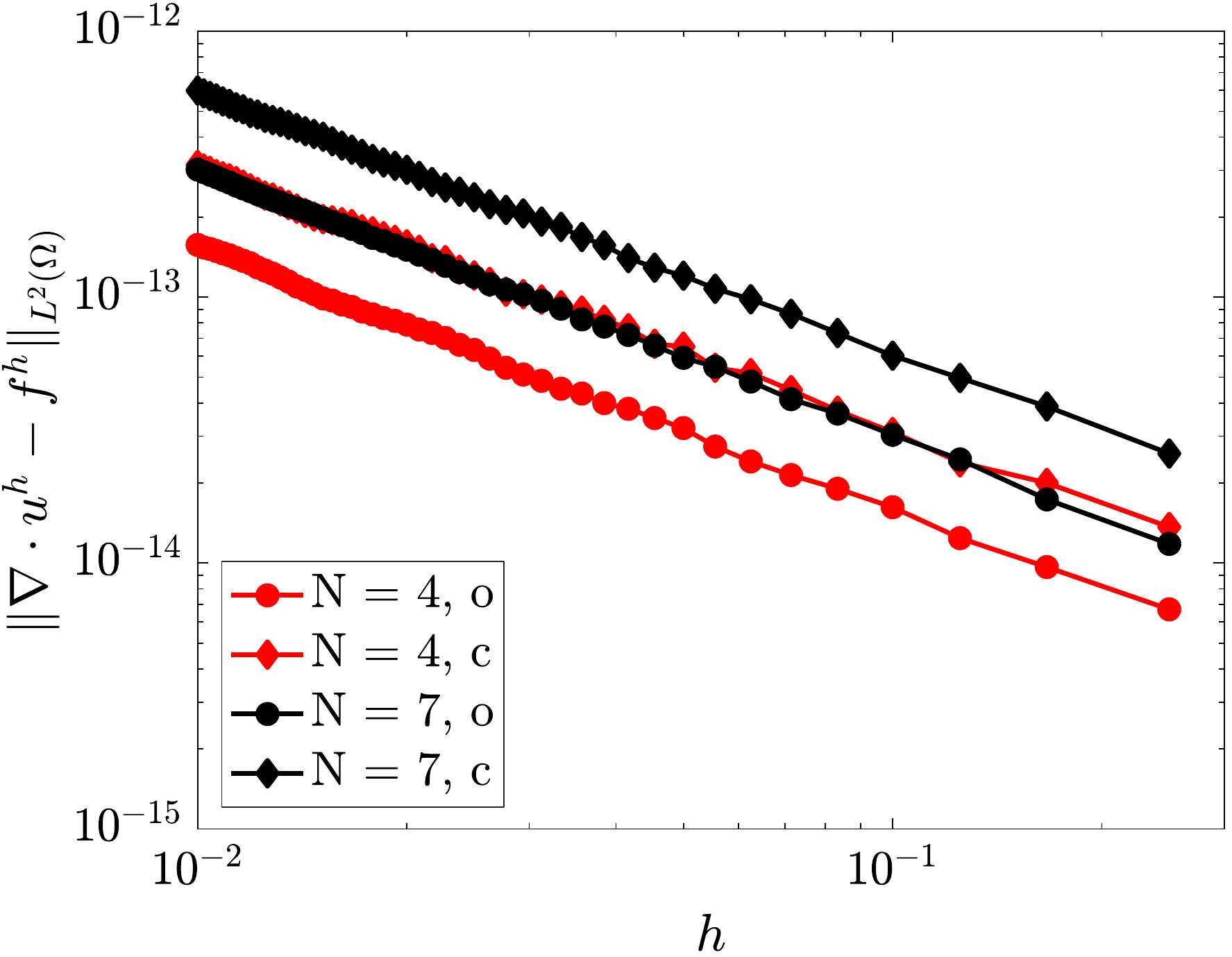}
\includegraphics[scale=0.3]{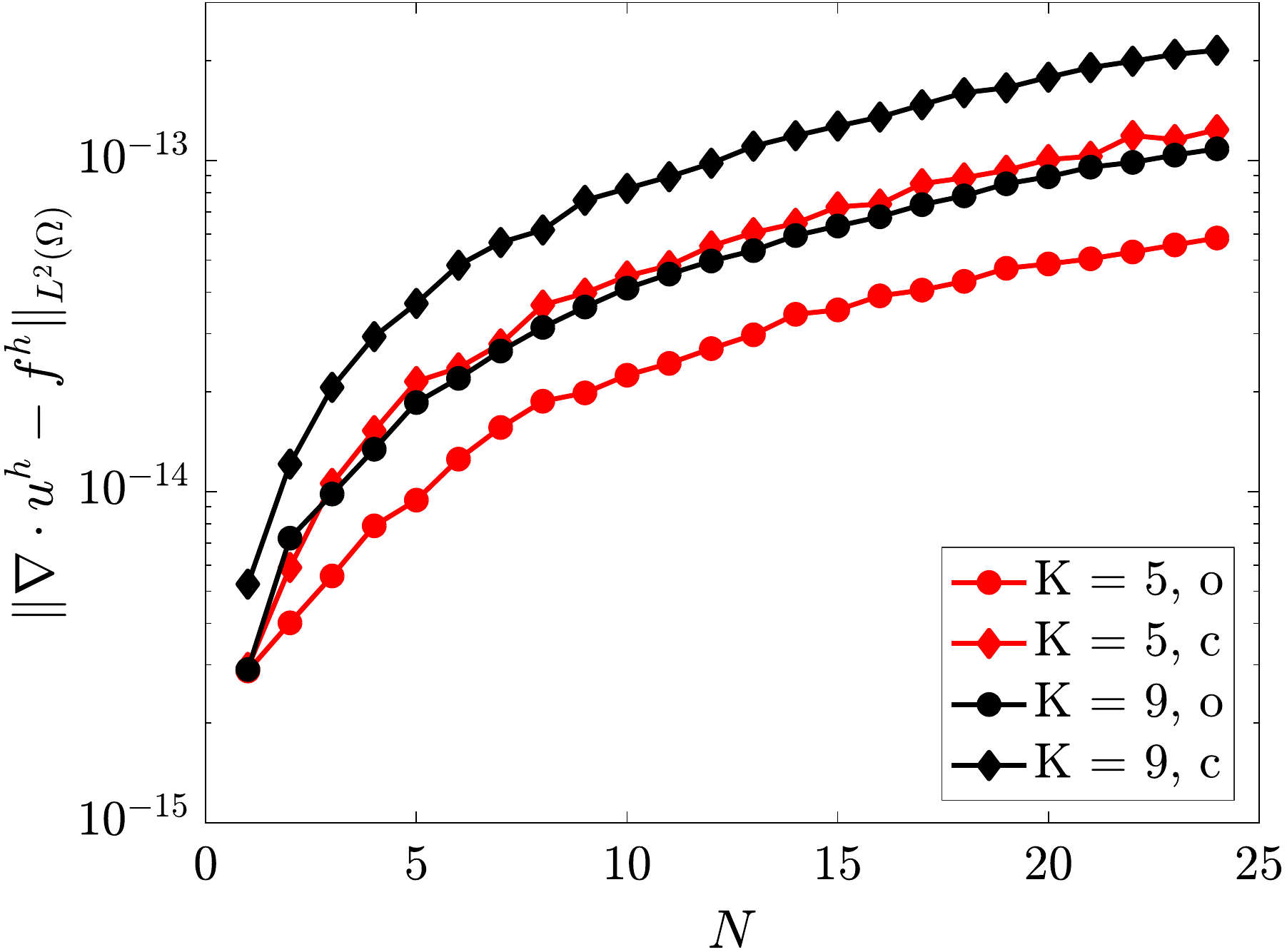}
\caption{$L^2$-error in divergence of velocity: Left: $h$-refinement;  Right: $N$-refinement.}
\label{fig:mass}
\end{figure}

In Fig. \ref{fig:mass} we show the $L^2$ - error for $\Vert \nabla \cdot u^h - f^h \Vert$.
On the left side as a function of element size, $h = 1 / \sqrt{K}$, and on the right side as a function of polynomial degree of the basis functions.
In both cases the maximum error observed is of $\mathcal{O}\bb{10^{-12}}$.
\begin{figure}
\includegraphics[scale=0.3]{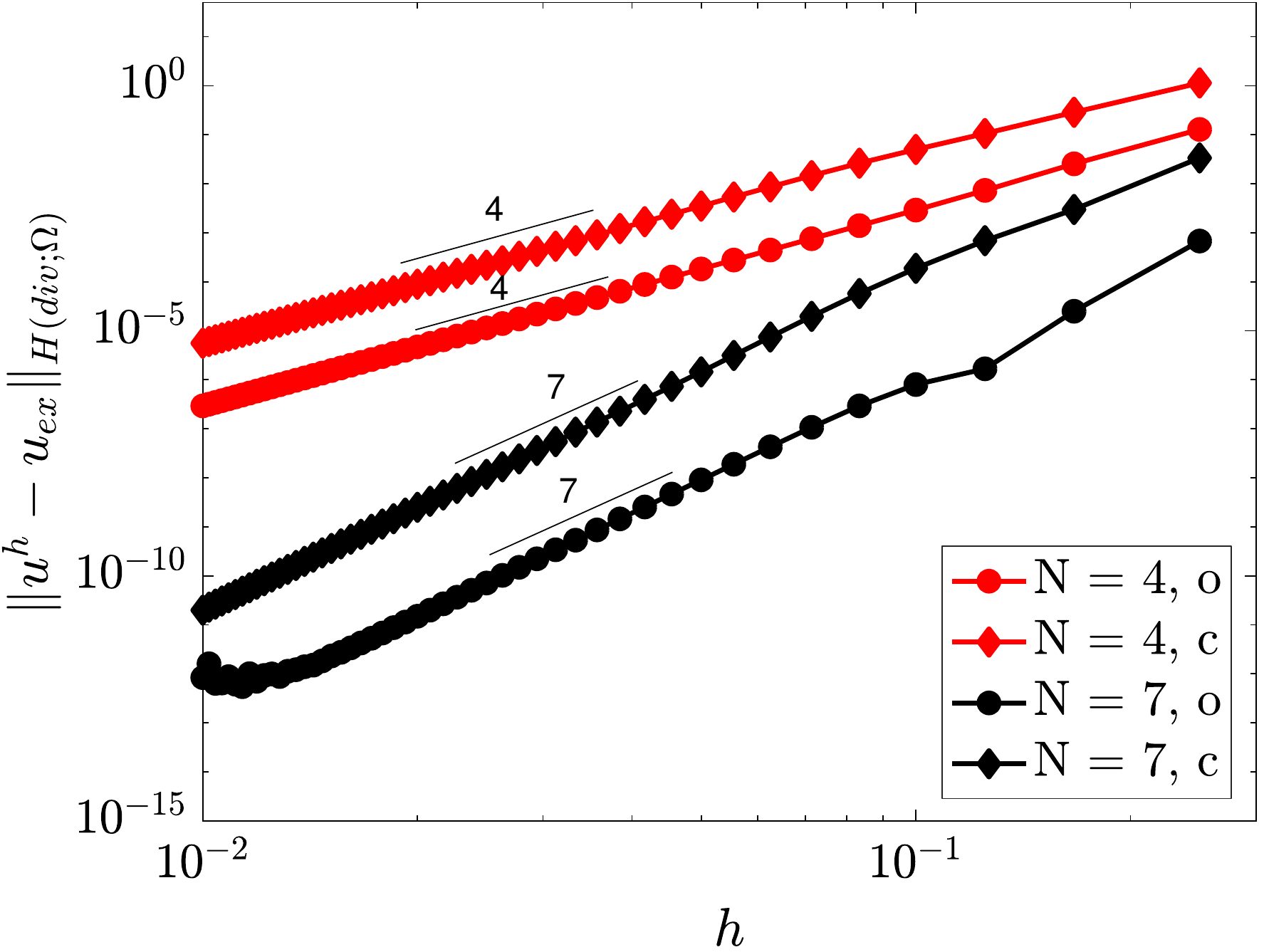}
\includegraphics[scale=0.3]{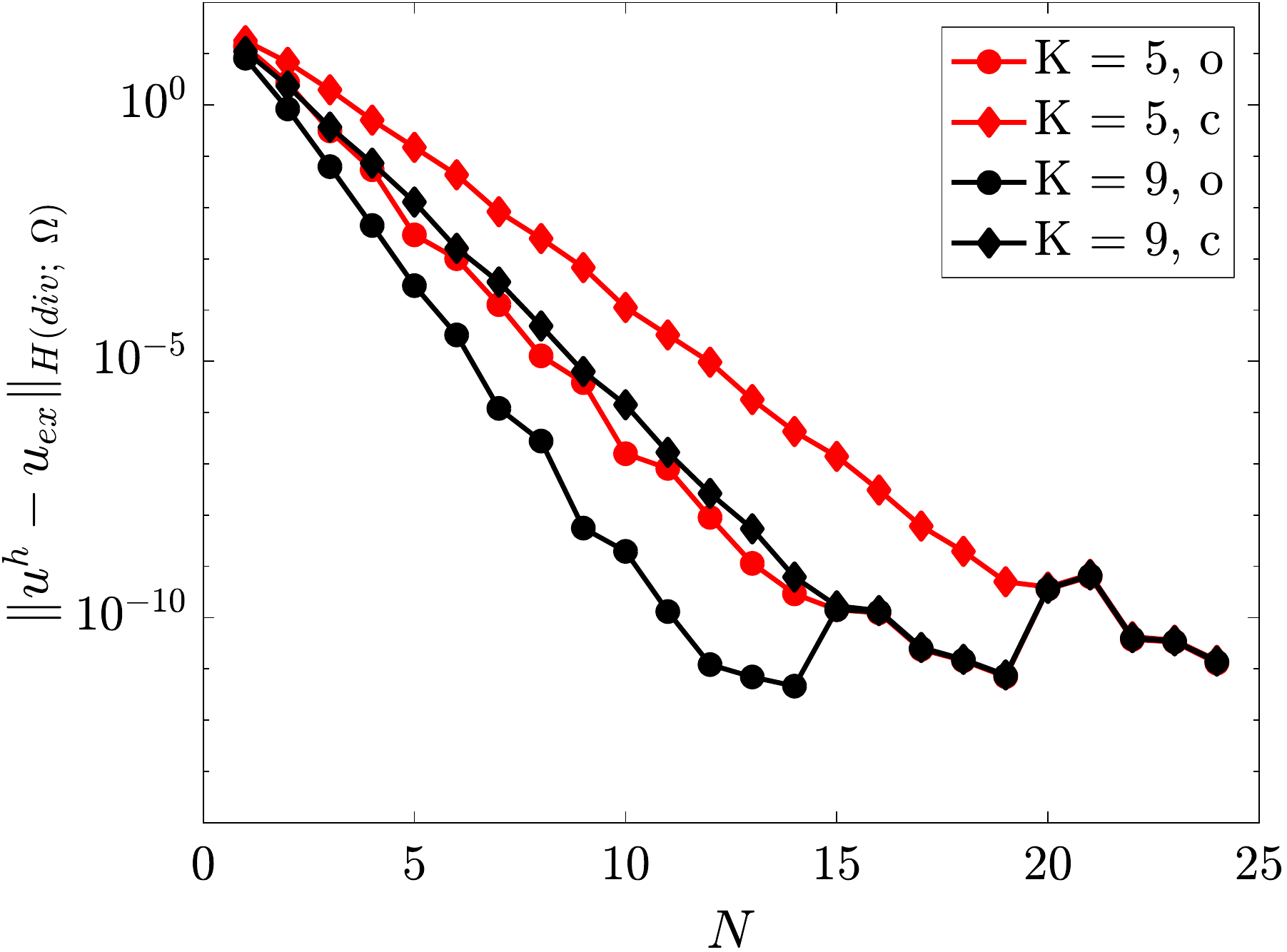} \nl
\includegraphics[scale=0.3]{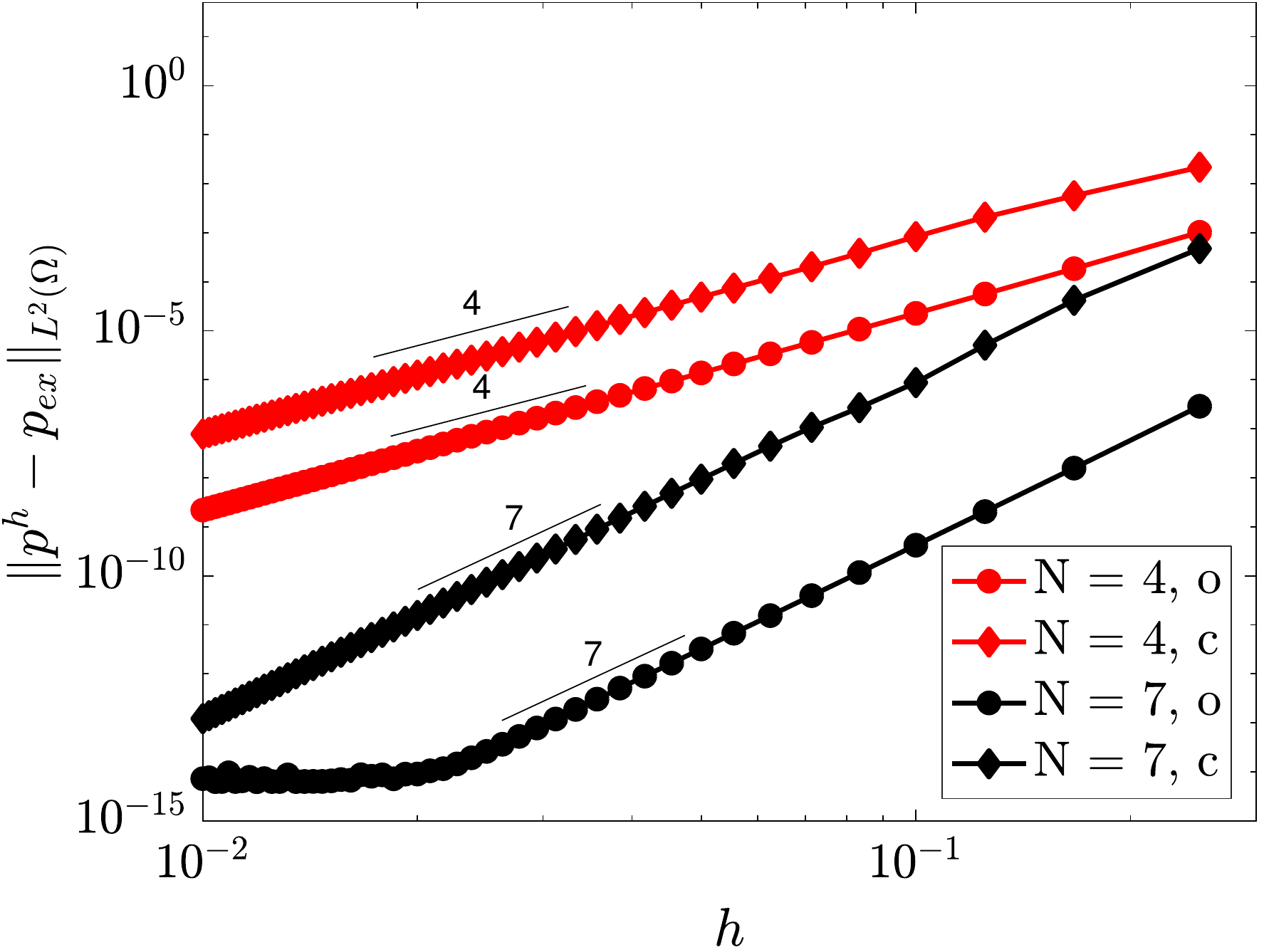}
\includegraphics[scale=0.3]{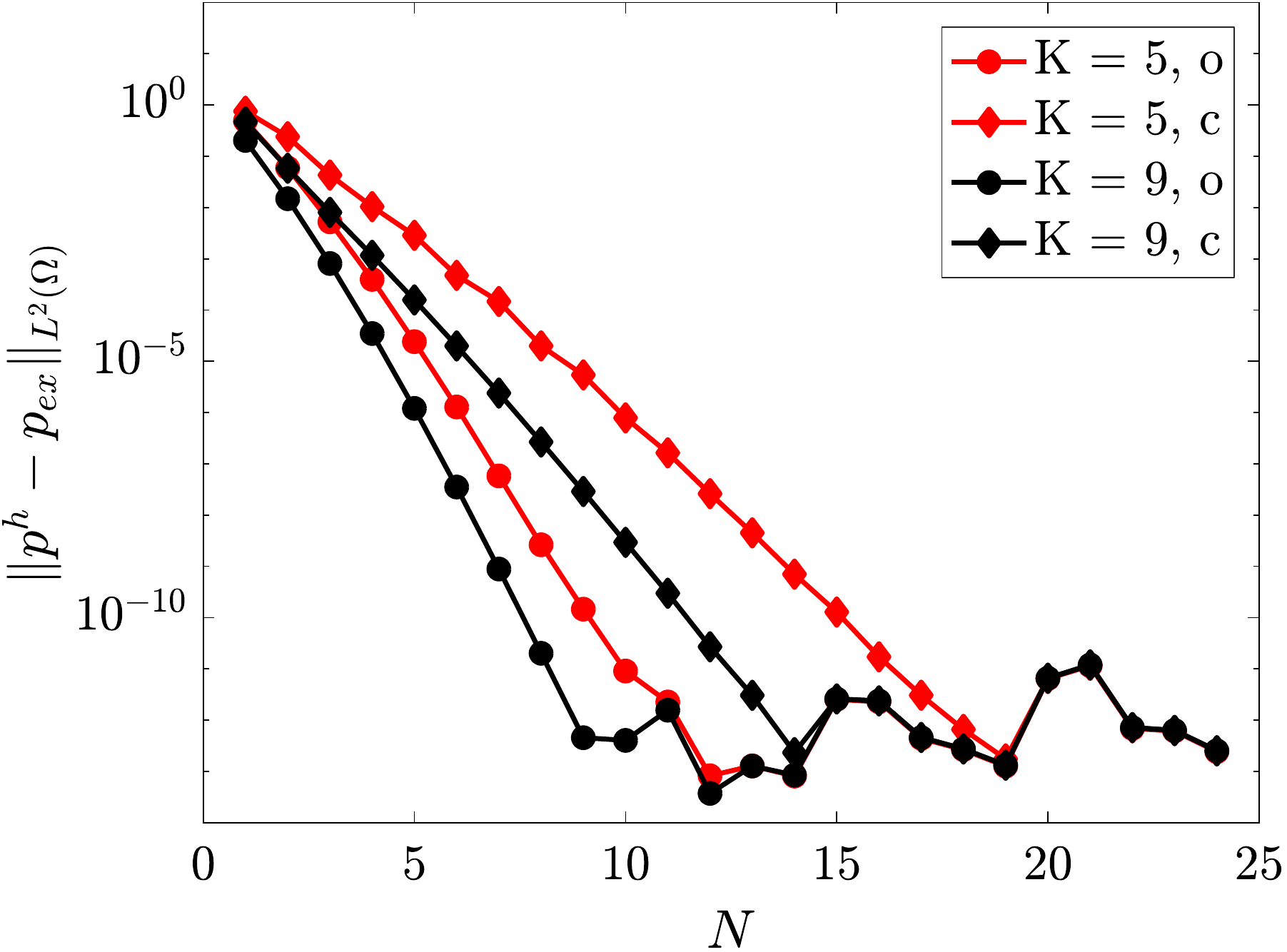}
\caption{Top row: Error in $H\bb{\mathrm{div}; \Omega}$ norm for velocity; Bottom row: $L^2$-error in pressure. Left: $h$-refinement;  Right: $N$-refinement.}
\label{fig:flux}
\end{figure}

In Fig. \ref{fig:flux}, on the top two figures we show the error in the $H\bb{\mathrm{div};\Omega}$ norm for the velocity; and at the bottom two figures we show the error in $L^2\bb{\Omega}$ norm for the pressure.
On the left we have $h$-convergence plots, and on the right we have $N$-convergence plots.
In all the figures, for the same number of elements, $K$, and polynomial degree, $N$, the error is higher for the curved mesh.

On the left we see that the error decreases with the element size.
The slope of error rate of convergence is $N$, which is optimal for both curved and orthogonal meshes.
On the right we see exponential convergence of the error with increasing polynomial degree of basis for both orthogonal and curved meshes.
\bibliographystyle{siam}
\bibliography{./library_darcy.bib}
\end{document}